# ASYMPTOTICS IN RESPONSE-ADAPTIVE DESIGNS GENERATED BY A TWO-COLOR, RANDOMLY REINFORCED URN

BY CATERINA MAY AND NANCY FLOURNOY

*Università degli Studi di Milano, Università del Piemonte Orientale and University of Missouri–Columbia*

This paper illustrates asymptotic properties for a response-adaptive design generated by a two-color, randomly reinforced urn model. The design considered is optimal in the sense that it assigns patients to the best treatment, with probability converging to one. An approach to show the joint asymptotic normality of the estimators of the mean responses to the treatments is provided in spite of the fact that allocation proportions converge to zero and one. Results on the rate of convergence of the number of patients assigned to each treatment are also obtained. Finally, we study the asymptotic behavior of a suitable test statistic.

**1. Introduction.** The present paper is devoted to studying asymptotic properties of sequential, response-adaptive designs generated by a two-color, generalized Pólya urn that is reinforced every time it is sampled with a random number of balls that are the same color as the ball that was extracted. The model is based on the two-color, randomly reinforced urn studied in [22] and [18], and for conciseness we will denote it by the acronym *RRU*. An *RRU* generalizes to discrete or continuous responses the urn models proposed initially in [9] and [17] for dichotomous responses, and which were also applied to select an optimal dosage in [8]. This work has been stimulated by the fact that a design driven by an *RRU* allocates units to the best response with a probability that converges to one. In the context of a response-adaptive design used to allocate patients in a clinical trial, this property is desirable from an ethical point of view; for this reason, the results obtained in this work will be illustrated within a clinical trial framework.

However, the reader should note that response-adaptive designs are of fundamental importance in many areas of applications, for instance in industrial problems.

The experimenter has two simultaneous goals: collecting evidence to determine the superior treatment, and biasing the allocations toward the better treatment in order to reduce the proportion of subjects in the experiment that receive the inferior treatment. Patients enter the experiment sequentially and are allocated randomly to a treatment, according to a rule that depends on the previous allocations and the previous observed responses. A vast number of adaptive designs have been proposed in recent years; informed reviews are found in [28] and in [13]. Many of them are based on generalized urn models; see, for instance, [27], which traces an historical development of generalized urn models, their properties and applications in sequential designs. Among more recent literature on urn models and adaptive designs, there are, for instance, [3, 4, 16, 30]. Notwithstanding their generality, their assumptions don't cover the $RRU$. A major reason is that the mean replacement matrix of an $RRU$ is not irreducible. Although in the past many response-adaptive designs have been focused on binary responses, more attention recently has been given to continuous outcomes. Among others, we note the procedures proposed in [2, 14, 29].

Consider a clinical trial conducted to compare two competing treatments, say $B$ and $W$, and a response-adaptive design. Indicate, with $N_B(n)$ and $N_W(n)$, the number of patients allocated through the $n$th patient to treatment $B$ and $W$, respectively. Many of these designs allocate patients targeting a certain proportion $\rho \in (0,1)$; the proportions of patients $N_B(n)/n$ and $N_W(n)/n$ allocated to each treatment converge almost surely to $\rho$ and $1-\rho$, respectively, where $\rho$ may be ad hoc or may be determined by some optimality criteria, which are usually a function of the unknown parameters of the outcomes. The adaptive design considered in this paper is different, because its optimality property is to assign patients to the best treatment with a proportion that converges almost surely to 1, while the proportion of patients allocated to the inferior treatment converges to 0. When the two treatments are equivalent, the design allocates the proportion of patients with a random limit in $[0,1]$.

After the specification of the model and the provision of preliminary results in Section 2, the first part of this work, included in Section 3, is dedicated to the study of the exact rate of convergence to infinity of the sample sizes $N_B(n)$ and $N_W(n)$. Because the treatment given corresponds to the color of the ball that is drawn, $N_B(n)$ and $N_W(n)$ also correspond to the number of balls of each color sampled from an $RRU$ through the $n$th draw. Moreover, we obtain the order of convergence of the process representing the proportion of black balls contained in the urn at every stage.

The asymptotic properties of response-adaptive designs studied in literature are usually based on the hypothesis that the target allocation $\rho$ is a determined value in $(0, 1)$. In some of those procedures, despite the randomness of the number of patients $N_B(n)$ and $N_W(n)$ allocated to treatments $B$ and $W$, respectively, and the complex dependence structure of the random variables involved, it has been proved that joint normality of the estimators of the mean responses based on the observed responses still holds. An important contribution is given in [21], where the authors provide a general method to prove consistency and an easy, general, nonmartingale approach to prove asymptotic normality of estimators based on adaptively observed allocations. As the authors show, their method can be applied to a wide class of adaptive designs targeting an allocation proportion $\rho \in (0, 1)$. Although their basic framework covers the adaptive design considered in this paper and strong consistency of the adaptive estimators of the parameters involved is derived from it, their method can't work for proving asymptotic normality when $\rho$ is exactly equal to one or $\rho$ may have a random behavior as in the $RRU$ procedure. Therefore, in Section 4, we prove that the joint asymptotic normality of the estimators of the mean responses still holds, both when the two treatments are equivalent and when one treatment is superior. The argument used resorts to a martingale technique that involves the concept of *mixing convergence*; mixing convergence is required to obtain the distribution of the two-sample t-statistic.

In Section 5, we consider the following hypothesis test: the experimenter wants to test the null hypothesis that the two treatments are equivalent, in the sense that the mean responses $\mu_B$ and $\mu_W$ are equal, against the alternative hypothesis that $\mu_B > \mu_W$. For this reason, we are interested in the distribution of the usual test statistic $\zeta_0$ for comparing the difference of the means, based on the observed responses, both under the null and the alternative hypothesis. When the target allocation is a value $\rho$ in $(0, 1)$, from the joint normality of the estimators of the means, it can be deduced using Slutsky's theorem that the test statistic $\zeta_0$ still has an asymptotic normal distribution. Also, under the alternative hypothesis, $\zeta_0$ has the same asymptotic noncentrality parameter as in the classical case, in which the sample sizes $n_B$ and $n_W$ are deterministic. For a review of the approach to these situations, see [12, 13, 29].

In the response-adaptive procedure considered in this paper, since the limit allocation is a random variable under the null hypothesis, we can't apply directly Slutsky's Theorem to derive the asymptotic normality of the test statistic $\zeta_0$ from the joint normality result of estimators. Notwithstanding, it is proved that under the null hypothesis the asymptotic normality of $\zeta_0$ holds in this procedure. The proof uses the mixing property of the convergence as described in Section 4. Moreover, under the alternative hypothesis, it is showed that $\zeta_0$ is a specific mixture of normal distributions, and we characterize its representation. A final discussion concludes the paper. To make it easy for the readers, proofs are located in the Appendix.

**2. Model specification and preliminary results.** Let $\{(Y_B(n), Y_W(n)): n \geq 1\}$ be a sequence of independent and identically distributed random response vectors with marginal distributions $\mathcal{L}_B$ and $\mathcal{L}_W$, discrete or continuous on $\mathbb{R}$. Consider an urn containing initially $b$ black balls and $w$ white balls, where $b$ and $w$ are two strictly positive real numbers. With the arrival of the first patient ($n = 1$), a ball is drawn at random from the urn and its color is observed. We define a random variable $\delta_1$ that we assume to be independent of the potential response vector $(Y_B(i), Y_W(i))$ for every $i \geq 1$ such that $\delta_1 = 1$ if the extracted ball is black, while $\delta_1 = 0$ if the extracted ball is white. So $\delta_1$ is a Bernoulli random variable with parameter $Z_0 = b/(b+w)$. After the ball is extracted, if it is black, it is replaced in the urn together with $U(Y_B(1))$ black balls. Otherwise, if it is white, it is replaced in the urn together with $U(Y_W(1))$ white balls, where $U$ is an arbitrary measurable function such that $U(Y_B(1))$ and $U(Y_W(1))$ have distribution on a nonnegative and bounded real set. In typical applications, $U$ will be a monotone function. (Note that $U$ may be chosen as the identity function when the distributions $\mathcal{L}_B$ and $\mathcal{L}_W$ have nonnegative and bounded support.)

This process is then reiterated at every instant $n+1$, $n \geq 1$. A ball is extracted and we define a random variable $\delta_{n+1}$, indicating its color $\delta_{n+1} = 1$ if the ball extracted is black and $\delta_{n+1} = 0$ if the ball extracted is white. We always assume that $\delta_{n+1}$ is independent of the potential response vector $(Y_B(i), Y_W(i))$ for every $i \geq n+1$. After the ball is extracted, it is replaced in the urn together with

$$\delta_{n+1} U(Y_B(n+1)) + (1 - \delta_{n+1}) U(Y_W(n+1))$$

balls of the same color. So, given the $\sigma$-algebra

$$(2.1) \quad \begin{aligned} \mathcal{F}_n = \sigma(\delta_1, \delta_1 Y_B(1) \\ + (1 - \delta_1) Y_W(1), \ldots, \delta_n, \delta_n Y_B(n) + (1 - \delta_n) Y_W(n)), \end{aligned}$$

$\delta_{n+1}$ is Bernoulli distributed with parameter

$$Z_n = \frac{B_n}{B_n + W_n},$$

where

$$\begin{cases} B_n = b + \sum_{i=1}^{n} \delta_i U(Y_B(i)), \\ W_n = w + \sum_{i=1}^{n} (1 - \delta_i) U(Y_W(i)). \end{cases}$$

The *RRU* procedure drives the allocations $\{\delta_n\}$. When $\delta_n$ is 1, allocate the $n$th patient to the first treatment, say treatment $B$, and let the random variable $Y_B(n)$ be the potential response of $n$th patient to treatment $B$. When $\delta_n$ is 0, allocate the $n$th patient to the second treatment, say treatment $W$, and let $Y_W(n)$ be the potential response of $n$th patient to treatment $W$. Only one response for the $n$th patient will be actually observed, which we can write as $Y(n) = \delta_n Y_B(n) + (1 - \delta_n) Y_W(n)$.

We thus generate the following processes: the sequence $\{\delta_n : n \geq 1\}$ of Bernoulli random variables and the sequence $\{Z_n : n \geq 0\}$ of random variables in $[0,1]$, representing the proportion of black balls present in the urn at every stage. Now the number of black balls and white balls that have been drawn from the urn through the $n$th treatment allocation can be written as $N_B(n) = \sum_{i=1}^{n} \delta_i$ and $N_W(n) = \sum_{i=1}^{n} (1 - \delta_i)$, respectively; clearly, $N_B(n) + N_W(n) = n$. Also note that $B_n$ and $W_n$ are the cumulative (transformed) observed responses to treatment $B$ and $W$, respectively, adjusted by the initial numbers of balls in the urn. We call $B_n$ and $W_n$ *cumulative responses* for short.

Let $m_B = \int U(y) \mathcal{L}_B(dy)$ and $m_W = \int U(y) \mathcal{L}_W(dy)$ be the means of the transformed responses. Then, from [5] and [22], we have the following limit for the proportion of black balls in the urn.

THEOREM 2.1. *If $m_B > m_W$, then $\lim_{n \to \infty} Z_n = 1$, almost surely.*

Suppose now that the sequences of responses $\{Y_B(n)\}$ and $\{Y_W(n)\}$ have finite means $\mu_B = \int y \mathcal{L}_B(dy)$ and $\mu_W = \int y \mathcal{L}_W(dy)$ and that, for instance, the treatment $B$ is preferred to the treatment $W$ if $\mu_B > \mu_W$. Then, choosing a function $U$ such that $\mu_B > \mu_W$ if and only if $m_B > m_W$ and $\mu_B = \mu_W$ if and only if $m_B = m_W$, Theorem 2.1 ensures that a *RRU*-design allocates patients to the superior treatment with probability converging to one as $n$ goes to infinity.

REMARK 1. As suggested in [23], the existence of a bounded function $U$ such that
$$\int y \mathcal{L}_B(dy) > \int y \mathcal{L}_W(dy) \Leftrightarrow \int U(y) \mathcal{L}_B(dy) > \int U(y) \mathcal{L}_W(dy)$$
and
$$\int y \mathcal{L}_B(dy) = \int y \mathcal{L}_W(dy) \Leftrightarrow \int U(y) \mathcal{L}_B(dy) = \int U(y) \mathcal{L}_W(dy)$$
is guaranteed by the theory of utility. In the situation illustrated in this paper, the experimenter expresses a preference among the distributions of responses in terms of the ordering of their means. Indeed, he may also express a different preference between $\mathcal{L}_B$ and $\mathcal{L}_W$. The theory of utility gives conditions which guarantee the existence of a bounded utility function $U$ such that the expected utilities of the elements of a class of probability distributions on $\mathbb{R}$ are ordered in the same way as a certain preference among the probability distributions (see, e.g., [7]).

2.1. *Preliminary results.* The process $\{Z_n : n \geq 0\}$ of the proportions of black balls is of primary interest for the study of stochastic processes generated by this particular generalization of Pólya urn. Moreover, in a *RRU*-design, $Z_n$ represents the conditional probability of allocating the $n$th patient to treatment $B$; the asymptotic behavior of this process is also essential for analyzing the asymptotic normality of estimators for these designs, as will be made clear in the next sections. [22] provides the following general result.

PROPOSITION 2.2. *The sequence of proportions $\{Z_n : n \geq 0\}$ is eventually a bounded super or sub-martingale. Therefore, it converges almost surely to a random limit $Z_\infty$ in $[0, 1]$.*

When the urn is reinforced by the random variables $U(Y_B(n))$ and $U(Y_W(n))$ with means $m_B$ and $m_W$ such that $m_B > m_W$, then, as given by Theorem 2.1, the limit $Z_\infty$ is equal to 1 almost surely.

Consider the case in which $\mathcal{L}_B = \mathcal{L}_W$ so that the urn reinforcements $U(Y_B(n))$ and $U(Y_W(n))$ have the same distribution, say $\mu$. Then, in [20], it is shown that $P(Z_\infty = x) = 0$ for every $x \in [0, 1]$. However, the exact distribution of $Z_\infty$ is unknown except in a few particular cases. When $\mu$ is a point mass at a non-negative real number $m$, the $RRU$ degenerates to Polya's urn, in which case $Z_\infty$ has a Beta$(b/m, w/m)$ distribution. This is also the case for binary responses (success/failure) when $m$ balls are added to the urn after each success is obtained. In fact, [1] proves that only the nonnull part of the reinforcement distribution needs to be considered. For the general $RRU$ with $\mathcal{L}_B = \mathcal{L}_W$, [1] characterizes the distribution of $Z_\infty$ as the unique continuous solution, satisfying some boundary conditions, of a specific functional equation in which the unknowns are distribution functions on $[0, 1]$.

When $m_B = m_W$, it may happen that $\int U(y)^k \mathcal{L}_B(dy) \neq \int U(y)^k \mathcal{L}_W(dy)$ for some $k \geq 2$, and then $U(Y_B(n))$ and $U(Y_W(n))$ may have different distributions. This is of particular interest because it corresponds to a situation in which the two treatments are considered equivalent. In Section 3 of this paper, we establish a fundamental property of $Z_\infty$ when $m_B = m_W$, that is, $P(Z_\infty = 1) = P(Z_\infty = 0) = 0$ in this case.

The following preliminary result regarding the limiting sample sizes on $B$ and $W$ is important for showing the asymptotic normality of the common statistic for testing differences in mean responses.

PROPOSITION 2.3. *$N_B(n)$ and $N_W(n)$ converge to infinity almost surely as $n \to \infty$.*

Now, let

$$\tau_n = \inf\left\{k : \sum_{i=1}^k \delta_i = n\right\} \quad \text{and} \quad \nu_n = \inf\left\{k : \sum_{i=1}^k (1-\delta_i) = n\right\}.$$

Thus, $\tau_n = j$ indicates that the $n$th observed response to treatment $B$ occurs for patient $j$, and $\{Y_B(\tau_n)\}$ is the subsequence of the potential responses $\{Y_B(n)\}$ that are the observed responses to $B$. Similarly, $\nu_n = j$ indicates that the $n$th observed response to treatment $W$ occurs for patient $j$, and the subsequence of potential responses $\{Y_W(n)\}$ to $W$ that are the observed responses is $\{Y_W(\nu_n)\}$. [6] proves independence properties of the sequences of observed responses, so that the strong consistency of estimators based on those sequences can be deduced.

PROPOSITION 2.4. *The sequences $\{Y_B(\tau_n)\}$ and $\{Y_W(\nu_n)\}$ are i.i.d. with distributions $\mathcal{L}_B$ and $\mathcal{L}_W$, respectively, and are independent one of each other.*

As a consequence of this proposition, we can model the observed responses of random sizes $N_B(n)$ and $N_W(n)$ to treatments $B$ and $W$ as samples from two i.i.d. populations generated by $\mathcal{L}_B$ and $\mathcal{L}_W$, respectively. Assume that $\mathcal{L}_B$ and $\mathcal{L}_W$ depend on unknown parameters $\theta_B$ and $\theta_W$. We have the following.

COROLLARY 1. *Suppose that $\tilde{\theta}_B$ and $\tilde{\theta}_W$ are estimators of $\theta_B$ and $\theta_B$ based on $n$-dimensional samples from two independent i.i.d. sequences generated by $\mathcal{L}_B$ and $\mathcal{L}_W$, respectively. Let $\hat{\theta}_B$ and $\hat{\theta}_W$ be the corresponding estimators computed on observed responses through time $n$ in the RRU-design [which have random sample sizes $N_B(n)$ and $N_W(n)$, resp.]. If $\tilde{\theta}_B$ and $\tilde{\theta}_W$ converge a.s. to $\theta_B$ and $\theta_W$, respectively, then also $\hat{\theta}_B$ and $\hat{\theta}_W$ converge a.s. to $\theta_B$ and $\theta_W$.*

**3. Rates of convergence of $N_B(n)$ and $N_W(n)$.** The aim of this section is to study the rate of convergence to infinity of the sample size sequences $N_B(n) = \sum_{i=1}^n \delta_i$ and $N_W(n) = \sum_{i=1}^n (1-\delta_i)$. We also obtain the rate of the convergence of the process $Z_n$. Results for dichotomous responses can be found in [8]. First, we have the following.

PROPOSITION 3.1.

(3.1) $\quad \lim_{n \to \infty} \dfrac{N_B(n)}{n} = Z_\infty \quad a.s. \quad \text{and} \quad \lim_{n \to \infty} \dfrac{N_W(n)}{n} = 1 - Z_\infty \quad a.s.,$

where $Z_\infty$ is the limit of the process $\{Z_n\}$ representing the proportion of black balls in the urn.

This result doesn't provide the rate of convergence to infinity of $N_B(n)$ and $N_W(n)$, because the limits $Z_\infty$ and $1 - Z_\infty$ may have a point mass at 0. The next theorem gives the exact rate of convergence to infinity of $N_B(n)$ and $N_W(n)$.

THEOREM 3.2. *If $m_B > m_W$, then:*

(i) $\lim_{n \to +\infty} \frac{N_B(n)}{n} = 1$, *a.s.;*
(ii) *there exist a random variable $\eta^2$ with $P(0 < \eta^2 < \infty) = 1$ such that*

$$\lim_{n \to +\infty} \frac{N_W(n)}{n^{m_W/m_B}} = \eta^2 \qquad a.s.$$

*If $m_B = m_W$, then $P(Z_\infty = 0) = P(Z_\infty = 1) = 0$. Hence, the limits in Proposition 3.1 say that the rate of convergence to infinity of $N_B(n)$ and $N_W(n)$ is $n$, almost surely, in this case.*

Part (i) of Theorem 3.2 follows immediately from Theorem 2.1 and Proposition 3.1, since when $m_B > m_W$, then $Z_\infty = 1$ almost surely and $N_B(n)/n$ converges to $Z_\infty$. The rest of the theorem is a consequence of an auxiliary result concerning the relative convergence rates of the cumulative responses to treatments $B$ and $W$, which holds both when $m_B = m_W$ and $m_B > m_W$.

THEOREM 3.3. $B_n/(W_n^{m_B/m_W})$ *converges almost surely to a random variable $\psi$ with support in $(0, \infty)$.*

Finally, we obtain a relevant consequence of Theorem 3.2: when $m_B > m_W$, the exact order of convergence of the proportion of white balls in the urn is the same as was obtained for the proportion of patients allocated to treatment $W$.

COROLLARY 2. *Let $\eta^2$ be the random variable given by Theorem 3.2. If $m_B > m_W$, then*

$$\lim_{n \to +\infty} \frac{1 - Z_n}{n^{m_W/m_B - 1}} = \frac{m_W}{m_B} \eta^2 \qquad a.s.$$

REMARK 2. There is a direct relationship among the random variable $Z_\infty$ and the limits in Theorem 3.2 and Theorem 3.3. When $m_B = m_W$, then $Z_\infty = \frac{\psi}{\psi + 1}$ and $1 - Z_\infty = \frac{1}{\psi + 1}$. When $m_B > m_W$, then $\eta^2 = \frac{1}{m_W}(\frac{m_B}{\psi})^{m_W/m_B}$.

REMARK 3. Theorem 3.3 has been inspired by some ideas contained in [25]. Actually, even when $m_B > m_W$, Theorem 3.3 doesn't follow from [25], Lemma 3.5, readapted to our case, which would assert that $\log W_n / \log B_n \to m_W/m_B$ as $n \to +\infty$, almost surely. As a counterexample, consider the case in which the rate of convergence of $B_n$ is $n$, while the rate of convergence of $W_n$ is $n^{m_W/m_B} \log n$. (This is also a counterexample to the fact that their true Theorem 2.2 is not a consequence of [25], Lemma 3.5.)

**4. Asymptotic normality.** Consider estimation of the means $\mu_B$ and $\mu_W$ of the responses to treatments. We define the following estimators based on the observed responses through patient $n$, with random sample sizes $N_B(n)$ and $N_W(n)$, respectively:

$$(4.1) \quad \hat{Y}_B(n) = \frac{\sum_{i=1}^n \delta_i Y_B(i)}{N_B(n)} \quad \text{and} \quad \hat{Y}_W(n) = \frac{\sum_{i=1}^n (1-\delta_i) Y_W(i)}{N_W(n)}.$$

Corollary 1 and the strong law of large numbers ensure that the (4.1) are strongly consistent. In this section, we show that these estimators, appropriately standardized, are jointly asymptotically normal, despite the randomness of $N_B(n)$ and $N_W(n)$, their dependence, and the fact that they don't satisfy the classical assumption, as also is required in [21], Theorem 3.2, that $N_B(n)/n$ and $N_W(n)/n$ converge in probability to a constant in $(0,1)$.

Before proceeding, we need to recall the concept of *mixing convergence*, which was introduced by [26] and provides a particularly elegant approach to martingale central limit theory (see [10] and [11]).

DEFINITION 4.1. Consider a sequences of random vectors $\mathbf{Y}_n = (Y_1(n), \ldots, Y_p(n))$ on a probability space $(\Omega, \mathcal{F}, P)$ that converges in distribution to $\mathbf{Y} = (Y_1, \ldots, Y_p)$. We say that the convergence is *mixing* if, for every point $\mathbf{y}$ of continuity for the distribution function of $\mathbf{Y}$ and for every event $E \in \mathcal{F}$,

$$\lim_{n \to \infty} P(Y_1(n) \leq y_1, \ldots, Y_p(n) \leq y_p, E) = P(Y_1 \leq y_1, \ldots, Y_p \leq y_p)P(E).$$

The mixing property of the convergence will be essential to our study of the asymptotic distribution of the test statistic in the following section. During this section and the next one, we will assume that the distributions of the responses $\mathcal{L}_B$ and $\mathcal{L}_W$ have finite variances, which we indicate by $\sigma_B{}^2$ and $\sigma_W{}^2$. The next theorem establishes the joint asymptotic normality of estimators (4.1).

THEOREM 4.1. *Either when $\mu_B = \mu_W$ or when $\mu_B > \mu_W$,*

$$\left( \frac{\sqrt{N_B(n)}}{\sigma_B}(\hat{Y}_B(n) - \mu_B), \frac{\sqrt{N_W(n)}}{\sigma_W}(\hat{Y}_W(n) - \mu_W) \right) \to^d \mathcal{N}(\mathbf{0}, \mathbf{I}) \quad \text{(mixing)}.$$

**5. Testing hypothesis.** In this section we consider the hypothesis test on the mean responses

$$H_0 : \mu_B = \mu_B \quad \text{versus} \quad H_1 : \mu_B > \mu_W.$$

Observe that, from Remark 1, these hypotheses are equivalent to $H_0 : m_B = m_B$ versus $H_1 : m_B > m_W$. We characterize the classical statistic for two samples $\{Y_B(n)\}$ and $\{Y_W(n)\}$, i.i.d. and with law $\mathcal{L}_B$ and $\mathcal{L}_W$, respectively, when applied to the response-adaptive design that motivates this paper.

First, consider the usual two-sample test statistic observed in a fixed design with sample sizes $n_B$ and $n_W$, respectively,

$$(5.1) \qquad \zeta_0 = \frac{\bar{Y}_B(n_B) - \bar{Y}_W(n_W)}{\sqrt{s_B^2/n_B + s_W^2/n_W}},$$

where $\bar{Y}_B(n_B)$ and $\bar{Y}_W(n_W)$ are the sample means, and $s_B^2$ and $s_W^2$ are consistent estimators of the variances. Suppose that, for some $\rho \in (0,1)$, $n_B/(n_B + n_W) \to \rho$, $n_W/(n_B + n_W) \to 1 - \rho$ as $n_B$ and $n_W \to \infty$; then, from the classical central limit theorem,

$$(5.2) \qquad \left( \frac{\sqrt{n_B}}{\sigma_B}(\bar{Y}_B(n_B) - \mu_B), \frac{\sqrt{n_W}}{\sigma_W}(\bar{Y}_W(n) - \mu_W) \right) \longrightarrow^d \mathcal{N}(\mathbf{0}, I),$$

we can deduce that $\zeta_0$ converges in distribution to a standard normal random variable under the null hypothesis, while it is asymptotically normal with noncentrality parameter

$$(5.3) \qquad \phi = \frac{\mu_B - \mu_W}{\sqrt{\sigma_B{}^2/n_B + \sigma_W^2/n_W}} \approx \sqrt{n} \frac{\mu_B - \mu_W}{\sqrt{\sigma_B{}^2/\rho + \sigma_W{}^2/(1-\rho)}}$$

under the alternative hypothesis.

Now, consider a response-adaptive procedure with random sample sizes such that $N_B(n)/n \to \rho$ and $N_W(n)/n \to (1 - \rho)$, where $\rho$ is a determined value in $(0,1)$, even if unknown a priori. If the result (5.2), with $N_B(n)$ and $N_W(n)$ replacing $n_B$ and $n_W$, still holds, it can be deduced, similarly to the classical case and using Slutsky's theorem, that the test statistic $\zeta_0$ [with $N_B(n)$ and $N_W(n)$ replacing $n_B$ and $n_W$] still has asymptotic normal distribution, and the same asymptotic noncentrality parameter under the alternative hypothesis. A discussion on the power of the test and the sample size calculation in this context can be found in [13].

Examine now the $RRU$-design considered in this paper. Under the null hypothesis, $N_B(n)/n$ and $N_W(n)/n$ converge to random limits $Z_\infty$ and $(1 - Z_\infty)$, respectively. So, we can't derive the asymptotic normality of the test statistic as in the classical case; we need the mixing property of the convergence proved in Theorem 4.1. Also, under the alternative hypothesis we find a different situation: $N_B(n)/n \to 1$ and $N_W(n)/n \to 0$, almost surely.

This implies that, asymptotically, $\zeta_0(n)$ carries information about the value of the mean of one only treatment. Notwithstanding this loss of balance, the use of the test statistic $\zeta_0(n)$ is still reasonable because the rates of $N_B(n)$ and $N_W(n)$ carry information about the difference between $\mu_B$ and $\mu_W$.

Let us define

$$\hat{\sigma}_B^2(n) = \sum_{i=1}^n \delta_i(Y_B(i) - \hat{Y}_B(n))^2 / N_B(n)$$

and

$$\hat{\sigma}_W^2(n) = \sum_{i=1}^n (1 - \delta_i)(Y_W(i) - \hat{Y}_B(n))^2 / N_W(n),$$

which are strong consistent estimators of the variances $\sigma_B^2$ and $\sigma_W^2$ from Corollary 1.

THEOREM 5.1. *The random process*

$$(5.4) \quad \zeta(n) = \frac{\hat{Y}_B(n) - \hat{Y}_W(n) - (\mu_B - \mu_W)}{\sqrt{\sigma_B^2/N_B(n) + \sigma_W^2/N_W(n)}} \to^d N(0,1) \qquad (mixing),$$

*both under the null and under the alternative hypothesis. Moreover, the same result holds with $\hat{\sigma}_B^2(n)$ and $\hat{\sigma}_W^2(n)$ in (5.4), instead of $\sigma_B^2$ and $\sigma_W^2$.*

Let us examine the distribution of the test statistic

$$\zeta_0(n) = \frac{\hat{Y}_B(n) - \hat{Y}_W(n)}{\sqrt{\hat{\sigma}_B^2(n)/N_B(n) + \hat{\sigma}_W^2(n)/N_W(n)}}.$$

When $H_0$ is true, $\zeta_0(n)$ is asymptotically normal from Theorem 5.1. So, one can then construct the following critical region with asymptotic level of significance $\alpha$:

$$C_\alpha = \{\zeta_0(n) > z_{1-\alpha}\}.$$

Let $\eta$ be the positive square root of $\eta^2$. The following corollary establishes that, when the alternative hypothesis is true, the test statistic $\zeta_0(n)$ is a mixture of normal distributions and characterizes its representation.

COROLLARY 3. *Under the alternative hypothesis, the conditional distribution of $\zeta_0(n)$, given the random variable $\eta^2$ defined in Theorem 3.2, is asymptotically normal with mean equal to $\sqrt{n^{m_W/m_B}} \eta \frac{\mu_B - \mu_W}{\sigma_W}$ and unit variance.*

REMARK 4. Simulations on the mixing distribution $\eta$ are provided in [19] and [15]. In particular, it is shown that $\eta$ is not a point mass, and this implies, by Corollary 3, that $\zeta(n)$ is not asymptotically normal under the alternative hypothesis.

**6. Discussion.** In Section 3, we have proved that the rate of convergence of the number of patients assigned to the worst treatment is determined by the ratio of $m_W$ and $m_B$: the smaller the value of this quantity, the more slowly $N_W(n)$ increases. However, for smaller values of $m_W/m_B$, the probability of assigning patients to the best treatment converges to one faster. On the other hand, the study of the noncentrality parameter in Section 5 shows that small values of $m_W/m_B$ cause a loss of the power of the test for treatment mean differences. This observation generates open research questions regarding the best choice for the function $U$ that determines the values of $m_B$ and $m_W$.

The theoretical distribution of the test statistic $\zeta_0$ for mean responses comparison has been studied in Section 5: $\zeta_0$ is asymptotically normal under the null hypothesis, while it is a specific mixture of normal distributions when the alternative hypothesis holds true. Hence, the law of the mixing distribution $\eta$ is of fundamental importance for calculating the exact power of the test. The distribution of $\eta$ is an open problem; some of its properties are investigated by simulation in [19] and [15]. Without knowing the exact distribution of $\eta$, the power of the test statistic must be evaluated by simulation.

Simulated results concerning the performance of an $RRU$-design are provided in [24] and [23] (for normal responses), and in [15] (for binary responses). The number of treatment failures to achieve a determined power is compared with a default design. Simulations show that an $RRU$-design is a viable alternative when the difference between the mean responses to treatments has moderate to large values. The $RRU$-design is also appealing with binary responses when success probabilities are small.

This work needs to be extended to the K-treatment problem. In [8] and [22] limiting results for the proportion of balls of each type generated by the RRU in the K-treatment problem has been obtained. But they did not obtain associated distribution theory for comparing treatment means and dealing with the associated multiple comparisons. We have laid a foundation for such results in this paper, but we believe the required extension will contain technical difficulties, and so we leave it now as an open problem.

In the present work, we have solved some asymptotic theory problems generated by the fact that an $RRU$-design has a very desirable property that can't be approached with the usual methods presented in literature. We wish that this study may offer a contribution to the development of research in response-adaptive, optimal designs.

APPENDIX: PROOFS

PROOF OF PROPOSITION 2.3. Set $\beta$ such that $P(U(Y_W(1)) \leq \beta) = 1$. Following [22], Theorem 2, if $\tau = \inf\{n \geq 1 : \delta_n = 0\}$, then, for $k \geq 1$,

$$P(\tau > k) \leq \frac{b}{b+w} \frac{b+\beta}{b+w+\beta} \cdots \frac{b+(k-1)\beta}{b+w+(k-1)\beta}$$

$$= \exp\left(\sum_{n=0}^{k-1} \log\left(\frac{b+n\beta}{b+w+n\beta}\right)\right).$$

Since $\sum_{n=0}^{\infty} \log[(b+n\beta)/(b+w+n\beta)] = -\infty$, it follows that $\lim_{k\to\infty} P(\tau > k) = 0$, and hence $P(\tau < \infty) = 1$. From the strong Markov property we obtain that

$$P(\delta_n = 0, \text{ i.o. }) = 1,$$

and then $N_W(n) \to \infty$, a.s. The proof for $N_B(n)$ is similar. $\square$

PROOF OF PROPOSITION 3.1. Since $\mathbf{E}(\delta_i | \mathcal{F}_{i-1}) = Z_{i-1}$, and (from Proposition 2.3) $\sum_{i=1}^{n} \delta_i \to \infty$, almost surely, it follows from Levy's extension of the Borel–Cantelli lemma that, almost surely,

$$\sum_{i=1}^{n} Z_{i-1} \to \infty \quad \text{and} \quad \frac{\sum_{i=1}^{n} \delta_i}{\sum_{i=1}^{n} Z_{i-1}} \to 1.$$

Because $\lim_{n\to\infty} Z_n = Z_\infty$, almost surely, $(\sum_{i=1}^{n} Z_{i-1})/n$ converges to $Z_\infty$ by Cesaro's lemma. Hence, $(\sum_{i=1}^{n} \delta_i)/n$ converges almost surely to $Z_\infty$. $\square$

In order to prove Theorem 3.3, we need the following two lemmas.

LEMMA A.1. *If $m_B > m_W$, then:*

(i) *the rate of convergence of $B_n$ is $n$, almost surely;*
(ii) *the rate of convergence of $W_n$ is greater then $n^\gamma$ for some real constant $\gamma > 0$, almost surely.*

*If $m_B = m_W$, then:*

(iii) *the rate of convergence of $B_n + W_n$ is $n$, almost surely;*
(iv) *the rates of convergence of both $B_n$ and $W_n$ are greater than $n^\gamma$, almost surely, for some real constant $\gamma > 0$.*

PROOF. (i) This is an immediate consequence of Corollary 1 and of the fact that when $m_B > m_W$ then $N_B(n)/n \to 1$ (using Proposition 3.1 and Theorem 2.1).

$$\lim_{n\to+\infty} \frac{B_n}{n} = \lim_{n\to+\infty} \frac{b}{n} + \frac{\sum_{i=1}^{n} \delta_i U(Y_B(i))}{n}$$

$$= \lim_{n\to+\infty} \frac{b}{n} + \frac{\sum_{i=1}^{n} \delta_i U(Y_B(i))}{N_B(n)} \frac{N_B}{n} = m_B \quad \text{a.s.}$$

(ii) Let's consider the conditional increments of the process $B_n/W_n^\kappa$, for some $\kappa > m_B/m_W$.

$$\mathbf{E}\left(\frac{B_{n+1}}{W_{n+1}^\kappa} - \frac{B_n}{W_n^\kappa}\bigg|\mathcal{F}_n\right)$$

$$= \frac{B_n}{B_n + W_n}\mathbf{E}\left(\frac{B_n + U(Y_B(n+1))}{W_n^\kappa} - \frac{B_n}{W_n^\kappa}\bigg|\mathcal{F}_n\right)$$

$$+ \frac{W_n}{B_n + W_n}\mathbf{E}\left(\frac{B_n}{(W_n + U(Y_W(n+1)))^\kappa} - \frac{B_n}{W_n^\kappa}\bigg|\mathcal{F}_n\right)$$

$$= \frac{B_n}{B_n + W_n}\frac{\mathbf{E}(U(Y_B))}{W_n^\kappa}$$

$$+ \frac{W_n B_n}{B_n + W_n}\mathbf{E}\left(\frac{1}{(W_n + U(Y_W(n+1)))^\kappa} - \frac{1}{W_n^\kappa}\bigg|\mathcal{F}_n\right).$$

By a Taylor expansion of the function $f(x) = 1/(a+x)^\kappa$ with $x = U(Y_W(n+1))$ and $a = W_n$, we can choose a constant $c$ such that whenever $W_n \geq 1$

$$\mathbf{E}\left(\frac{1}{(W_n + U(Y_W(n+1)))^\kappa}\bigg|\mathcal{F}_n\right) \leq \frac{1}{W_n^\kappa} - \frac{\kappa}{W_n^{\kappa+1}}\left(\mathbf{E}(U(Y_W)) - \frac{c}{W_n}\right).$$

Thus, we obtain that

$$\mathbf{E}\left(\frac{B_{n+1}}{W_{n+1}^\kappa} - \frac{B_n}{W_n^\kappa}\bigg|\mathcal{F}_n\right)$$

(A.1)

$$\leq \frac{B_n}{B_n + W_n}\frac{1}{W_n^\kappa}\left(\mathbf{E}(U(Y_B)) - \kappa\mathbf{E}(U(Y_W)) + \frac{\kappa c}{W_n}\right).$$

From inequality (A.1), note that if $\kappa > \mathbf{E}(U(Y_B))/\mathbf{E}(U(Y_W))$, then the process $B_n/W_n^\kappa$ is (eventually) a positive supermartingale, and then it converges, almost surely, to a finite limit. Since from part (i) of the lemma $B_n/n \to m_B$ a.s., it follows that also $n/W_n^\kappa$ converges almost surely to a finite limit. Hence, for every $\varepsilon > 0$, $n/W_n^{\kappa+\varepsilon}$ converges a.s. to 0. This means that $W_n^{\kappa+\varepsilon} > n$ eventually, that is, $W_n > n^{1/(\kappa+\varepsilon)}$ a.s., eventually.

(iii) If $m_B = m_W$, then $\lim_{n \to +\infty}\frac{B_n + W_n}{n} = m_B$, almost surely, because

$$\frac{B_n + W_n}{n} = \lim_{n \to +\infty}\frac{b+w}{n} + \frac{\sum_{i=1}^n \delta_i U(Y_B(i))}{n}$$

$$+ \frac{\sum_{i=1}^n (1-\delta_i) U(Y_W(i))}{n}$$

$$= \lim_{n \to +\infty}\frac{b+w}{n} + \frac{\sum_{i=1}^n \delta_i U(Y_B(i))}{N_B(n)}\frac{N_B(n)}{n}$$

$$+ \frac{\sum_{i=1}^n (1-\delta_i) U(Y_W(i))}{N_W(n)}\frac{N_W(n)}{n},$$

which, by Corollary 1 and Proposition 3.1, converges almost surely to $m_B Z_\infty + m_W(1 - Z_\infty)$. Since $m_B = m_W$, this is equal to $m_B$.

(iv) From the part (iii) of the lemma, it follows that, eventually, $B_n + W_n > n \cdot m_B \cdot 3/4$ on a set of probability one. Let $F_B := \{B_n > n \cdot m_B/4, \text{eventually}\}$ and $F_W := \{W_n > n \cdot m_B/4, \text{eventually}\}$. By Proposition 2.2, $Z_n$ converges a.s. in $[0,1]$. Then, the ratio $B_n/W_n$ converges a.s. in $[0,+\infty]$ and, therefore, $P(F_B \cup F_W) = 1$. In fact, on $(F_B \cup F_W)^c$, we have that $\liminf_n B_n/n \le m_B/4$, $\liminf_n W_n/n \le m_B/4$ together with $B_n + W_n > n \cdot m_B \cdot 3/4$ that implies $\liminf_n B_n/W_n \le 1/2 < 2 \le \limsup_n B_n/W_n$, that is, $Z_n$ does not converges on $(F_B \cup F_W)^c$. Hence, $P(F_B \cup F_W)^c = 0$. What remains to prove is that, on $F_B$ (resp. on $F_W$), $W_n$ (resp. $B_n$) is eventually greater than $n^\gamma$, for some real constant $\gamma > 0$.

Let us focus on what happens on the set $F_B$. Using the same argument as in the proof of part (ii) of this lemma, it follows from equation (A.1) that for $\kappa > 1$ the process $B_n/W_n^\kappa$ is (eventually) a positive supermartingale, and then it converges almost surely to a finite limit. This implies that, for every $\varepsilon > 0$, $B_n/W_n^{\kappa+\varepsilon}$ converges to 0 and then

$$\text{(A.2)} \qquad B_n/W_n^{\kappa+\varepsilon} < 1 \qquad \text{on } F_B, \text{ eventually.}$$

Since, on $F_B$, $B_n > n \cdot m_B/4$ eventually, then, by (A.2), we deduce that, on $F_B$, $n/W_n^{\kappa+\varepsilon} \cdot m_B/4 < 1$, a.s., that is, $W_n > n^{1/(\kappa+\varepsilon)} \cdot (m_B/4)^{1/(\kappa+\varepsilon)}$, eventually. □

The second lemma is a general fact about convergence of random sequences; for lack of a better reference, see [25], Lemma 3.2.

LEMMA A.2. *Let $\{\xi_n : n \ge 0\}$ be a random sequence that is measurable with respect to a filtration $\{\mathcal{F}_n\}$. Define*

$$\Delta_n = \mathbf{E}(\xi_{n+1} - \xi_n | \mathcal{F}_n);$$
$$Q_n = \mathbf{E}((\xi_{n+1} - \xi_n)^2 | \mathcal{F}_n).$$

*If $\sum_n \Delta_n < +\infty$ and $\sum_n Q_n < +\infty$ on a set of probability one, then $\xi_n$ converges to a finite random variable, almost surely, as $n$ goes to infinity.*

PROOF OF THEOREM 3.3. We apply Lemma A.2 to the process

$$\xi_n = \log \frac{B_n}{W_n^{m_B/m_W}}$$

in order to prove that it converges almost surely to a finite random variable. This implies that $B_n/(W_n^{m_B/m_W})$ converges almost surely to a strictly positive and finite random variable.

$$\begin{aligned}\Delta_n &= \mathbf{E}(\xi_{n+1} - \xi_n | \mathcal{F}_n) \\ &= \mathbf{E}(\log B_{n+1} - \log B_n | \mathcal{F}_n) - \frac{m_B}{m_W}\mathbf{E}(\log W_{n+1} - \log W_n | \mathcal{F}_n) \\ &= \frac{B_n}{B_n + W_n}\mathbf{E}(\log(B_n + U(Y_B(n+1))) - \log B_n | \mathcal{F}_n) \\ &\quad - \frac{m_B}{m_W}\frac{W_n}{B_n + W_n}\mathbf{E}(\log(W_n + U(Y_W(n+1))) - \log W_n | \mathcal{F}_n) \\ &= \frac{B_n}{B_n + W_n}\mathbf{E}\left(\int_0^{U(Y_B(n+1))} \frac{1}{B_n + t}\,dt \Big| \mathcal{F}_n\right) \\ &\quad - \frac{m_B}{m_W}\frac{W_n}{B_n + W_n}\mathbf{E}\left(\int_0^{U(Y_W(n+1))} \frac{1}{W_n + t}\,dt \Big| \mathcal{F}_n\right).\end{aligned}$$

By a Taylor expansion of the function $f(x) = 1/(x+t)$, it follows that, for $B_n$ sufficiently large, there exist constants $c_1$ and $c_2$ such that, for every $t$

$$(\text{A.3}) \qquad \frac{1}{B_n} - c_2\frac{t}{B_n^2} \leq \frac{1}{B_n + t} \leq \frac{1}{B_n} - \frac{t}{B_n^2} + c_1\frac{t^2}{B_n^3},$$

and, similarly, for $1/(W_n + t)$. Hence, it follows that

$$(\text{A.4}) \qquad -\frac{1}{B_n + W_n}\left(\frac{h_1}{B_n} + \frac{h_2}{W_n^2}\right) \leq \Delta_n \leq \frac{1}{B_n + W_n}\left(\frac{k_1}{B_n^2} + \frac{k_2}{W_n}\right)$$

for some constant $k_1, k_2, h_1, h_2$. Thanks to the rates of convergence of $B_n$ and $W_n$ shown in Lemma A.1, we obtain that $\sum_n \Delta_n < +\infty$, a.s.

$$\begin{aligned}Q_n &= \mathbf{E}((\xi_{n+1} - \xi_n)^2 | \mathcal{F}_n) \\ &= \frac{B_n}{B_n + W_n}\mathbf{E}((\log(B_n + U(Y_B(n+1))) - \log B_n)^2 | \mathcal{F}_n) \\ &\quad + \frac{W_n}{B_n + W_n}\mathbf{E}((\log(W_n + U(Y_W(n+1))) - \log W_n)^2 | \mathcal{F}_n) \\ &= \frac{B_n}{B_n + W_n}\mathbf{E}\left(\left(\int_0^{U(Y_B(n+1))} \frac{1}{B_n + t}\,dt\right)^2 \Big| \mathcal{F}_n\right) \\ &\quad + \left(\frac{m_B}{m_W}\right)^2 \frac{W_n}{B_n + W_n}\mathbf{E}\left(\left(\int_0^{U(Y_W(n+1))} \frac{1}{W_n + t}\,dt\right)^2 \Big| \mathcal{F}_n\right).\end{aligned}$$

Since, for positive $t$, $\frac{1}{B_n + t} \leq \frac{1}{B_n}$ and $\frac{1}{W_n + t} \leq \frac{1}{W_n}$, we obtain

$$Q_n \leq \mathbf{E}\,\frac{B_n}{B_n + W_n}\left(\frac{E(U(Y_B)^2)}{B_n{}^2}\right) + \left(\frac{m_B}{m_W}\right)^2 \frac{W_n}{B_n + W_n}\left(\frac{E(U(Y_W)^2)}{W_n{}^2}\right);$$

hence, $\sum_n Q_n < +\infty$, a.s. $\square$

PROOF OF THEOREM 3.2. Let $m_B > m_W$. Observe that

$$\lim_{n \to +\infty} \frac{B_n}{W_n^{m_B/m_W}} = \lim_{n \to +\infty} \frac{b + \sum_{i=1}^n \delta_i U(Y_B(i))}{(w + \sum_{i=1}^n (1-\delta_i) U(Y_W(i)))^{m_B/m_W}}$$

$$= \lim_{n \to +\infty} \frac{\sum_{i=1}^n \delta_i U(Y_B(i))/N_B(n)}{(\sum_{i=1}^n (1-\delta_i) U(Y_W(i))/N_W(n))^{m_B/m_W}}$$

$$\times \frac{N_B(n)}{N_W(n)^{m_B/m_W}}.$$

From Theorem 3.3, we know that this limit is equal to $\psi \in (0, \infty)$. From Corollary 1, we have that $\sum_{i=1}^n \delta_i U(Y_B(i))/N_B(n)$ and $\sum_{i=1}^n (1-\delta_i) U(Y_W(i))/N_W(n)$ converge a.s. to $m_B$ and $m_W$, respectively, and $N_B(n)/n$ converges a.s. to 1. So, it follows that $N_W(n)^{m_B/m_W}/n$ converges a.s. to a random variable in $(0, \infty)$. Hence, we obtain part (ii) of the theorem.

Let $m_B = m_W$. From Proposition 3.1, we know that, almost surely, $N_B(n)/n \to Z_\infty$ and $N_W(n)/n \to 1 - Z_\infty$. Since, from Theorem 3.3, $B_n/W_n$ converges a.s. to $\psi \in (0, \infty)$, it follows that

$$Z_n = \frac{B_n}{B_n + W_n} = \frac{B_n/W_n}{B_n/W_n + 1}$$

converges a.s. to $\psi/(\psi + 1)$, which is then in $(0, 1)$. Hence,

(A.5) $$P(Z_\infty = 0) = P(Z_\infty = 1) = 0. \qquad \square$$

PROOF OF COROLLARY 2.

$$\frac{1 - Z_n}{n^{m_W/m_B - 1}} = \frac{1}{n^{m_W/m_B - 1}} \frac{W_n}{B_n + W_n}$$

$$= \frac{1}{n^{m_W/m_B - 1}} \frac{w + \sum_{i=1}^n (1 - \delta_i) U(Y_W(i))}{b + w + \sum_{i=1}^n \delta_i U(Y_B(i)) + \sum_{i=1}^n (1 - \delta_i) U(Y_W(i))}$$

$$= \frac{1}{n^{m_W/m_B - 1}}$$

$$\times \left( w + \left( \sum_{i=1}^n (1 - \delta_i) U(Y_W(i))/N_W(n) \right) N_W(n) \right)$$

$$\times \left( b + w + \left( \sum_{i=1}^n \delta_i U(Y_B(i))/N_B(n) \right) N_B(n) \right.$$

$$\left. + \left( \sum_{i=1}^n (1 - \delta_i) U(Y_W(i))/N_W(n) \right) N_W(n) \right)^{-1}$$

$$= \left( w/n^{m_W/m_B} + \sum_{i=1}^{n}(1-\delta_i)U(Y_W(i))/N_W(n) \cdot N_W(n)/n^{m_W/m_B} \right)$$

$$\times \left( (b+w)/n + \sum_{i=1}^{n}\delta_i U(Y_B(i))/N_B(n) \cdot N_B(n)/n \right.$$

$$\left. + \sum_{i=1}^{n}(1-\delta_i)U(Y_W(i))/N_W(n) \cdot N_W(n)/n \right)^{-1}.$$

Now, using Corollary 1 and Theorem 3.2, the numerator of the last equality converges almost surely to $\eta^2 m_W$, while the denominator converges almost surely to $m_B$. □

In order to prove Theorem 4.1, we need the following lemma, which extends [10], Corollary 3.2, to the multi-dimensional case. It can be obtained using the Cramér–Wold device as in [11], Theorem 12.6.

LEMMA A.3. *Let* $\{\mathbf{X}_{ni} = (X_{ni}{}^1, \ldots, X_{ni}{}^p), \mathcal{F}_{ni}, 1 \leq i \leq k_n\}$ *be a zero-mean, square-integrable, p-dimensional martingale array differences. Suppose that*

$$\sum_i \mathbf{E}[\|\mathbf{X}_{ni}\|^2 I(\|\mathbf{X}_{ni}\| \geq \epsilon)|\mathcal{F}_{n,i-1}] \to^P 0 \quad \text{for all } \epsilon > 0,$$

*and*

$$\mathbf{\Sigma}_n := \sum_i \mathbf{E}[(\mathbf{X}_{ni})'\mathbf{X}_{ni}|\mathcal{F}_{n,i-1}] \to^P \mathbf{\Sigma}.$$

*If* $\mathbf{\Sigma}$ *is a positive definite matrix with probability* 1, *then*

$$\sum_i \mathbf{X}_{ni} \mathbf{\Sigma}_n^{-1/2} \to^d \mathcal{N}(\mathbf{0}, \mathbf{I}) \quad (mixing).$$

PROOF OF THEOREM 4.1. Let $k_n = n$ and $\mathcal{F}_{ni} = \sigma(\delta_1, \delta_1 Y_B(1) + (1-\delta_1)Y_W(1), \ldots, \delta_i, \delta_i Y_B(i) + (1-\delta_i)Y_W(i), \delta_{i+1})$.

Since $\mathbf{E}(\delta_i(Y_B(i)-\mu_B)|\mathcal{F}_{n,i-1}) = \delta_i \cdot \mathbf{E}(Y_B(i)-\mu_B) = 0$ and $\mathbf{E}((1-\delta_i)(Y_W(i)-\mu_W)|\mathcal{F}_{n,i-1}) = (1-\delta_i) \cdot \mathbf{E}(Y_W(i)-\mu_W) = 0$, it follows that $\{\mathbf{X}_{ni}, \mathcal{F}_{ni}, 1 \leq i \leq k_n\}$, with

$$\mathbf{X}_{ni} = \left( \frac{\delta_i(Y_B(i)-\mu_B)}{\sqrt{n}}, \frac{(1-\delta_i)(Y_W(i)-\mu_W)}{\sqrt{n^{m_W/m_B}}} \right),$$

is a 2-dimensional martingale array differences. $\mathbf{X}_{ni}$ satisfies

$$\sum_i \mathbf{E}[\|\mathbf{X}_{ni}\|^2 I(\|\mathbf{X}_{ni}\| \geq \epsilon)|\mathcal{F}_{n,i-1}]$$

$$\leq \sum_i \mathbf{E}\left(\frac{\delta_i(Y_B(i) - \mu_B)^2}{n} I(|Y_B(i) - \mu_B| \geq \epsilon\sqrt{n})|\mathcal{F}_{n,i-1}\right)$$

$$+ \sum_i \mathbf{E}\bigg(\frac{(1-\delta_i)(Y_W(i) - \mu_W)^2}{n^{m_W/m_B}}$$

$$\times I(|Y_W(i) - \mu_W| \geq \epsilon\sqrt{n^{m_W/m_B}})|\mathcal{F}_{n,i-1}\bigg)$$

$$= \frac{N_B(n)}{n}\mathbf{E}((Y_B(1) - \mu_B)^2 I(|Y_B(1) - \mu_B| \geq \epsilon\sqrt{n}))$$

$$+ \frac{N_W(n)}{n^{m_W/m_B}}\mathbf{E}((Y_W(1) - \mu_W)^2 I(|Y_W(1) - \mu_W| \geq \epsilon\sqrt{n^{m_W/m_B}})) \to^P 0.$$

Moreover, it can be verified that

$$\mathbf{\Sigma}_n = \sum_i \mathbf{E}[(\mathbf{X}_{ni})'\mathbf{X}_{ni}|\mathcal{F}_{n,i-1}] = \begin{pmatrix} \frac{N_B(n)\sigma_B^2}{n} & 0 \\ 0 & \frac{N_W(n)\sigma_W^2}{n^{m_W/m_B}} \end{pmatrix}$$

and, by Theorem 3.2, $\mathbf{\Sigma}_n$ converges almost surely to a positive definite matrix $\mathbf{\Sigma}$ with probability 1, either when $m_B = m_W$ and when $m_B > m_W$, that is, when $\mu_B = \mu_W$ and $\mu_B > \mu_W$. We can then apply Lemma A.3 to $\mathbf{X}_{ni}$, and we obtain the theorem. $\square$

PROOF OF THEOREM 5.1. Let

$$\lambda_n = \frac{\sigma_B^2 N_W(n)/n}{\sigma_B^2 N_W(n)/n + \sigma_W^2 N_B(n)/n}.$$

Then, by Proposition 3.1,

$$\lambda_n \to \lambda_\infty := \frac{\sigma_B^2(1 - Z_\infty)}{\sigma_B^2(1 - Z_\infty) + \sigma_W^2 Z_\infty}$$

almost surely, and

$$\zeta(n) = \sqrt{\lambda_n} \cdot \frac{\sqrt{N_B(n)}}{\sigma_B}(\hat{Y}_B(n) - \mu_B) + \sqrt{1 - \lambda_n} \cdot \frac{\sqrt{N_W(n)}}{\sigma_W}(\hat{Y}_W(n) - \mu_W)$$

$$= \sqrt{\lambda_\infty} \cdot \frac{\sqrt{N_B(n)}}{\sigma_B}(\hat{Y}_B(n) - \mu_B) + \sqrt{1 - \lambda_\infty} \cdot \frac{\sqrt{N_W(n)}}{\sigma_W}(\hat{Y}_W(n) - \mu_W)$$

$$+ (\sqrt{\lambda_n} - \sqrt{\lambda_\infty}) \cdot \frac{\sqrt{N_B(n)}}{\sigma_B}(\hat{Y}_B(n) - \mu_B)$$

$$+ (\sqrt{1 - \lambda_n} - \sqrt{1 - \lambda_\infty}) \cdot \frac{\sqrt{N_W(n)}}{\sigma_W}(\hat{Y}_W(n) - \mu_W)$$

$$= \sqrt{\lambda_\infty} \cdot \frac{\sqrt{N_B(n)}}{\sigma_B}(\hat{Y}_B(n) - \mu_B) + \sqrt{1-\lambda_\infty} \cdot \frac{\sqrt{N_W(n)}}{\sigma_W}(\hat{Y}_W(n) - \mu_W)$$
$$+ o(1).$$

According to Theorem 4.1, it follows that $\zeta(n)$ converges mixing in distribution to a random variable having a characteristic function as

$$\mathbf{E}(\exp(-\tfrac{1}{2}\lambda_\infty t^2) \cdot \exp(-\tfrac{1}{2}(1-\lambda_\infty)t^2)) = e^{-t^2/2},$$

which is the characteristic function of a standard normal random variable. For the case of unknown variances, define

$$\gamma(n) = \sqrt{\frac{\sigma_B^2 N_W(n)/n + \sigma_W^2 N_B(n)/n}{\hat{\sigma}_B^2(n) N_W(n)/n + \hat{\sigma}_W^2(n) N_B(n)/n}},$$

which converges almost surely to 1, by Theorem 3.2. By noting that

$$\frac{\hat{Y}_B(n) - \hat{Y}_W(n) - (\mu_B - \mu_W)}{\sqrt{\hat{\sigma}_B^2(n)/N_B(n) + \hat{\sigma}_W^2(n)/N_W(n)}} = \zeta(n) \cdot \gamma(n),$$

from the mixing convergence of $\zeta(n)$ and Slutsky's theorem, we obtain the thesis. $\square$

PROOF OF COROLLARY 3. When $\mu_B > \mu_W$, we have that $\zeta_0(n) = \zeta(n) + \phi(n)$, where

$$\phi(n) = \frac{\mu_B - \mu_W}{\sqrt{\hat{\sigma}_B^2(n)/N_B(n) + \hat{\sigma}_W^2(n)/N_W(n)}}$$

is the noncentrality parameter. Since, for Theorem 3.2, $N_B(n)/n^{m_W/m_B} \to \infty$ and $N_W(n)/n^{m_W/m_B} \to \eta^2$ almost surely, then, for $n$ large, we can approximate the noncentrality parameter with

$$\phi(n) \approx \sqrt{n^{m_W/m_B}} \sqrt{\eta^2} \frac{\mu_B - \mu_W}{\sigma_W}.$$

For the mixing property of the convergence in Theorem 5.1, $\zeta(n)$ is asymptotically independent of $\eta^2$, and then we obtain the thesis. $\square$

**Acknowledgments.** Special thanks to P. Secchi, G. Aletti, J. Moler and C. Heyde for many discussions and suggestions that have been an essential contribution to this work. We thank also an anonymous referee for his very constructive comments to improve Sections 4 and 5.

Dipartimento di Matematica "F. Enrigues"
Università degli Studi di Milano
Via Cesare Saldini 50
20133 Milan
Italy
and
Dipartimento di Scienze Econimiche
  e Metodi Quantitativi
Università del Piemonte Orientale
Via Perrone 18
28100 Novara
Italy
E-mail: cmay@mat.unimi.it

Department of Statistics
University of Missouri–Columbia
Columbia, Missouri 65211
USA
E-mail: flournoyn@missouri.edu